\input amstex
\magnification=\magstep1 
\baselineskip=13pt
\documentstyle{amsppt}
\vsize=8.7truein \CenteredTagsOnSplits \NoRunningHeads
\def\den{\operatorname{den}}
\def\UU{\Cal U}
\def\DD{\Cal D}
\def\Pr{\bold{P}}
\def\EE{\bold{E}}
\topmatter
 
 \title Searching for dense subsets in a graph via the partition function  \endtitle 
\author Alexander Barvinok and Anthony Della Pella \endauthor
\address Department of Mathematics, University of Michigan, Ann Arbor,
MI 48109-1043, USA \endaddress
\email barvinok$\@$umich.edu, adellape$\@$umich.edu \endemail
\date July 2018 \enddate
\thanks  The first author was partially supported by NSF Grant DMS 1361541.
\endthanks 
\abstract  For a set $S$ of vertices of a graph $G$, we define its density $0 \leq \sigma(S) \leq 1$ as the ratio of the number of edges of $G$ spanned by the vertices of $S$ to ${|S| \choose 2}$. We show that, given a graph $G$ with $n$ vertices and an integer $m$, the partition function $\sum_S \exp\{ \gamma m \sigma(S) \}$, where the sum is taken over all 
$m$-subsets $S$ of vertices and $0 < \gamma <1$ is fixed in advance,  can be approximated within relative error $0 < \epsilon < 1$ in quasi-polynomial $n^{O(\ln m - \ln \epsilon)}$ time.
We discuss numerical experiments and observe that for the random graph $G(n, 1/2)$ one can afford a much larger $\gamma$, provided the ratio $n/m$ is sufficiently large.
\endabstract
\endtopmatter
\document

\head 1. Introduction and main results \endhead 

Let $G=(V, E)$ be an undirected graph, without loops or multiple edges. For a non-empty subset $S \subset V$ of vertices, we define the {\it density} $\sigma(S)$ as the fraction of the pairs of vertices of $S$ that span an edge of $G$: 
$$\sigma(S)={\left| {S \choose 2} \cap E\right| \over {|S| \choose 2}},$$
where ${S \choose 2}$ is the set of all unordered pairs of vertices from $S$. Hence $0 \leq \sigma(S) \leq 1$ for all subsets, $\sigma(S)=0$ if $S$ is an {\it independent set} and 
$\sigma(S)=1$ if $S$ is a {\it clique}. 

We are interested in the following general problem: given a graph $G=(V, E)$ with $|V|=n$ vertices and an integer $m \leq n$, estimate the highest density of an 
$m$-subset $S \subset V$. This is, of course, a hard problem: for example, testing whether a given graph contains a clique of a given size, or even estimating the size of the largest clique within a factor of $n^{1-\epsilon}$ for any $\epsilon >0$, fixed in advance, is already an NP-hard problem \cite{Ha99}, \cite{Zu99}. Moreover, modulo some plausible complexity assumptions, it is 
hard to approximate the highest density of an $m$-subset for a given $m$, within a constant factor, fixed in advance \cite{Bh12}. The best known efficient approximation 
achieves the factor of $n^{1/4}$ in quasi-polynomial $n^{O(\ln n)}$ time \cite{B+10}. There are indications that the factor $n^{1/4}$ might be hard to beat \cite{B+12}.
We note that the most interesting case is when $m$ grows and $n \gg m$, since the highest density of an $m$-subset can be computed in polynomial time up to an additive error of $\epsilon n^2/m^2$ for any $\epsilon >0$, fixed in advance \cite{FK99} (and if $m$ is fixed in advance, the densest $m$-subset can be found by the exhaustive search in polynomial time).

\subhead (1.1) Partition function \endsubhead In this paper, we approach the problem of finding the densest, or just a reasonably dense subset, via computing the 
{\it partition function}
$$\den_m(G; \gamma)={n \choose m}^{-1} \sum\Sb S \subset V: \\ |S| =m \endSb \exp\left\{ \gamma m \sigma(S) \right\},  \tag1.1.1$$
where $\gamma > 0$ is a parameter. The {\it exponential tilting}, $\sigma(S) \longmapsto \exp\left\{\gamma m \sigma(S) \right\}$,  see for example, Section 13.7 of \cite{Te99}, puts greater emphasis on the sets of higher density. Let us consider the set ${V \choose m}$ of all $m$-subsets of $V$ as a probability space with the uniform measure. It is not hard to see that for any $0 < \sigma_0 < 1$, we have  
$$ \sigma_0 + { \ln  \Pr\left(\sigma(S) \ \geq \ \sigma_0 \right) \over \gamma m} \ \leq \  { \ln \den_m(G; \gamma) \over \gamma m} \ \leq \ \max\Sb S \subset V: \\ |S|=m \endSb \sigma(S), \tag1.1.2$$
so the larger $\gamma$ we can afford, the better approximation for the densest $m$-subset we get. In particular, if we could choose $\gamma \gg \ln n$ then 
from (1.1.2) we could approximate the highest density of an $m$-subset within an arbitrarily small additive error. 

The partition function (1.1.1) was introduced in \cite{Ba15}, where an 
algorithm of quasi-polynomial $n^{O(\ln m -\ln \epsilon)}$ complexity was constructed to compute (1.1.1) within relative error $0<\epsilon<1$, when $\gamma=0.07$ and when $\gamma=0.27$, under additional assumptions that $n \geq 8m$ and $m \geq 10$. It follows from (1.1.2) that if the probability to hit an $m$-subset $S$ of density at least $\sigma_0$ at random is $e^{-o(m)}$ then 
we can certify the existence of an $m$-subset of density at least $\sigma_0-o(1)$ in quasi-polynomial time, just by computing (1.1.1). It is also shown in \cite{Ba15} that by successive conditioning, one can find in quasi-polynomial time an $m$-subset $S$ with density at least as high as certified by the value of (1.1.1).

In this paper, we present an algorithm, which, for any $0 < \gamma < 1$, fixed in advance, and a given $0 < \epsilon <1$, computes the value of (1.1.1) within relative error $\epsilon$ in quasi-polynomial $n^{O(\ln m -\ln \epsilon)}$ time, provided 
$n > \omega(\gamma) m$ for some constant $\omega(\gamma)>1$. This improvement from $\gamma=0.27$ to an arbitrary $\gamma <1$ required the addition of some new ideas to the technique of \cite{Ba15}. We still don't know, however, if (1.1.1) can be efficiently computed for {\it any} $\gamma >0$, fixed in advance, and as we remarked above, it is unlikely that (1.1.1) can be efficiently computed for $\gamma \gg \ln n$. Our numerical experiments seem to indicate that we can afford a substantially larger $\gamma$. This can be partially explained by the fact that for the Erd\H{o}s-R\'enyi random graph $G(n, 0.5)$ indeed a much larger $\gamma$ can be used with high probability.

\subhead (1.2) Multivariate partition function \endsubhead Given $n \times n$ symmetric complex matrix $Z=\left(z_{ij}\right)$ and $2 \leq m \leq n$, we define 
$$P_m(Z)=\sum\Sb S \subset \{1, \ldots, n\} \\ |S|=m \endSb \exp\left\{ \sum\Sb \{i, j\} \subset S \\ i \ne j \endSb z_{ij} \right\}. \tag1.2.1$$
Note that the diagonal entries of $Z$ are irrelevant, so we assume that $z_{ii}=0$ for all $i$.

Given a graph $G=(V, E)$ with set $V=\{1, \ldots, n\}$ of vertices and $\gamma>0$, we define $Z_0=\left(z_{ij}\right)$ by 
$$z_{ij}=\cases {\gamma \over m-1} &\text{if\ } \{i, j\} \in E \\ -{\gamma \over m-1} &\text{if\ } \{i, j\} \notin E, \endcases$$
and observe that 
$$\split P_m(Z_0)= &\sum\Sb S \subset \{1, \ldots, n\} \\ |S|=m \endSb \exp\left\{ m \gamma \sigma(S) - {\gamma m \over 2} \right\}\\=
&\exp\left\{ -{\gamma m \over 2}\right\}{n \choose m} \den_m(G; \gamma). \endsplit \tag1.2.2$$
Hence to compute (1.1.1) it suffices to compute $P_m(Z_0)$. We compute $P_m(Z_0)$ by interpolation, see \cite{Ba15}, \cite{Ba16}. For that, it suffices to show that $P_m(Z) \ne 0$ 
in some neighborhood of a path connecting the zero matrix to $Z_0$ in the space of complex matrices. 

We prove the following result.
\proclaim{(1.3) Theorem} For any $0 < \delta < 1$ there exist $\eta=\eta(\delta) >0$ and $\omega=\omega(\delta) >1$ such that if 
$n \geq \omega m$ then $P_m(Z) \ne 0$ for any $n \times n$  symmetric complex matrix $Z=\left(z_{ij}\right)$ such that 
$$\left| \Re\thinspace z_{ij}\right| \ \leq \ {\delta \over m-1} \quad \text{and} \quad \left| \Im\thinspace z_{ij} \right| \ \leq \ {\eta \over m-1} \quad \text{for all} \quad 
1 \leq i \ne j \leq n.$$
\endproclaim

We prove Theorem 1.3 in Sections 2 and 3. Using Theorem 1.3, in Section 4 we  present an algorithm of quasi-polynomial $n^{O(\ln m)}$ complexity to compute 
$P_m(Z_0)$ and hence $\den_m(G; \gamma)$ for any $0 < \gamma <1$, fixed in advance.

In \cite{Ba15} it was established that $P_m(Z) \ne 0$ in a polydisc 
$$\DD_{m,n}=\left\{ Z=\left(z_{ij}\right): \ |z_{ij}| \leq {0.27 \over m-1} \quad \text{for all} \quad 1 \leq i \ne j \leq n \right\}$$
provided $n \gg m$ and $m$ is large enough.
In Theorem 1.3, we establish that $P_m(Z) \ne 0$ in a more ``economical" domain, ``stretched" along the real part of the complex space of matrices.
 This allows us to improve the constant $\gamma$ for which 
$\den_m(G; \gamma)$ is still efficiently computable.

In Section 5, we discuss some results of our numerical experiments, which seem to indicate that we can afford an essentially bigger $\delta$ in Theorem 1.3. This can be partially explained by the fact that for the Erd\H{o}s-R\'enyi random graph $G(n, 0.5)$ this is indeed the case. Namely, we prove the following result in Section 6.

\proclaim{(1.4) Theorem} Let us choose positive integers $n$ and $2 \leq m \leq n$. For $n \times n$ symmetric matrix $W=\left(w_{ij}\right)$ of independent random variables, 
where 
$$\Pr\left(w_{ij}=1\right)=\Pr\left(w_{ij}=-1\right)={1 \over 2},$$ we define the polynomial
$$h_W(z)={n \choose m}^{-1}\sum\Sb S \subset \{1, \ldots, n\} \\ |S|=m \endSb \prod_{\{i, j\} \subset S} \left(1+ z w_{ij}\right).$$
Let $r >0$ and $\tau >1$ be real numbers. If 
$n \ \geq \ 2m^2 \left(1+r^2\right)^m +2m$
then the probability that $h_W(z)$ has a root in the disc $|z| \ < \ r/ \sqrt{2 \tau}$
does not exceed $1/\tau$.
\endproclaim
In particular, if $n \gg m^2$ then with high probability $h_W(z)$ has no roots in the disc 
$|z| < c/\sqrt{m}$,
for an arbitrary large $c>0$, fixed in advance. Similarly, if $\ln n \gg m$ then with high probability $h_W(z)$ has no roots in the disc 
$|z| < c$
for an arbitrary large $c >0$, fixed in advance. 

The polynomial $h_W(z)$ is easily translated into the partition function $\den_m(G; \gamma)$, where $G$ is the graph with set $V=\{1, \ldots, n\}$ of vertices and two vertices $\{i, j\}$ span an edge if and only if $w_{ij}=1$: for $0 < \alpha < 1$, we have 
$$h_W(\alpha)=(1-\alpha)^{m \choose 2} \den_m(G; \gamma) \quad \text{where} \quad \gamma={m-1 \over 2} \ln {1+ \alpha \over 1-\alpha}.$$
Consequently, with high probability we can can compute $\den_m(G; \gamma)$ in quasi-polynomial time for $\gamma$ as large as $\gamma=\sqrt{m}$ provided $n \gg m^2$ and 
as large as $\gamma=m$ provided $\ln n \gg m$. Since the graphs we experimented on were to a large degree random (but not necessarily Erd\H{o}s-R\'enyi $G(n, 0.5)$), we may have obtained overly optimistic numerical evidence.

\head 2. Preliminaries \endhead

We consider the partition function $P_m$  of Section 1.2 within a family of partition functions, which will allow us to prove Theorem 1.3 by induction.
\subhead (2.1) Functionals $P_{\Omega}(Z)$ \endsubhead Let us fix integers $n$ and $2 \leq m \leq n$. For a subset $\Omega \subset \{1, \ldots, n\}$ and $n \times n$ complex symmetric matrix $Z=\left(z_{ij}\right)$, we define 
$$P_{\Omega}(Z)=\sum \Sb S \subset \{1, \ldots, n\}: \\ |S|=m, \Omega \subset S \endSb \exp\left\{ \sum\Sb \{i, j\} \subset S \\ i\ne j \endSb z_{ij} \right\}$$
where we agree that $P_{\Omega}(Z) =0$ if $|\Omega| > m$.
In other words, we restrict the sum (1.2.1) defining $P_m(Z)$ onto subsets $S$ containing a given set $\Omega$. In particular, 
$$P_{\Omega}(Z)=P_m(Z) \quad \text{if} \quad \Omega=\emptyset.$$
The induction will be built on the following straightforward formulas:
$$P_{\Omega}(Z)= {1 \over m-|\Omega|} \sum_{j \in \{1, \ldots, n\} \setminus \Omega} P_{\Omega \cup \{j\}}(Z) \quad \text{provided} \quad |\Omega| < m \tag2.1.1$$
and for $i \ne j$, we have
$${\partial \over \partial z_{ij}} P_{\Omega}(Z)= \cases P_{\Omega}(Z) &\text{if\ } i, j \in \Omega,\\ P_{\Omega \cup \{j\}}(Z) &\text{if\ } i \in \Omega, j \notin \Omega, \\
P_{\Omega \cup \{i\}}(Z) &\text{if\ } i \notin \Omega, j \in \Omega, \\ P_{\Omega \cup\{i, j\}} (Z) &\text{if\ } i, j \notin \Omega. \endcases \tag2.1.2$$

We will often consider complex numbers as vectors in the plane, by identifying ${\Bbb C}={\Bbb R}^2$ and measuring, in particular, angles between non-zero complex numbers. 
We will use the following geometric lemma.

\proclaim{(2.2) Lemma}  Let $u_1, \ldots, u_n \in {\Bbb C}$ be non-zero complex numbers such that the angle between any two does not exceed $\theta$ for 
some $0 < \theta < \pi/2$. Suppose that
$$\Im \left(\sum_{j=1}^n u_j \right) =0 \quad \text{and} \quad \sum_{j=1}^n \left| u_j \right| =a.$$
Then 
$$\sum_{j=1}^n \left| \Im\thinspace u_j \right| \ \leq \ a\sin {\theta \over 2}.$$
\endproclaim
\demo{Proof} Scaling $u_j$, if necessary, without loss of generality we assume that $a=1$.

Without loss of generality, we assume that $\arg u_j \ne 0$ for $j=1, \ldots, n$. Indeed, if $\arg u_j =0$ for some $j$, we can remove the vector from 
the collection, which would make the sum 
$$\sum_{j=1}^n |u_j| \tag2.2.1$$
only smaller. Rescaling $u_j \longmapsto \tau u_j$ for some real $\tau >1$, we make (2.2.1) equal to 1 and increase 
$$\sum_{j=1}^n \left| \Im \thinspace u_j \right|. \tag2.2.2$$
Reflecting the vectors $u_j$ in the coordinate axis if necessary, without loss of generality we may assume that $\Re\thinspace u_1 \geq 0$ and $\Im\thinspace u_1 >0$.
Hence there is a vector, say $u_2$, such that $\Im\thinspace u_2 < 0$. We necessarily have $\Re\thinspace u_2 \geq 0$, since otherwise the angle between $u_1$ and $u_2$ exceeds $\pi/2$. Then for any vector $u_j$, we must have $\Re\thinspace u_j \geq 0$, since otherwise one of the angles formed by $u_j$ with $u_1$ or $u_2$ will 
exceed $\pi/2$. 

Hence without loss of generality, we assume that $\Re\thinspace u_j >0$ for $j=1, \ldots, n$. Let 
$$\alpha= \max_{j=1, \ldots, n} \arg u_j,$$
so that 
$$0 \ < \ \alpha \ < \ \theta$$
and let 
$$-\beta= \min_{j=1, \ldots, n} \arg u_j \ < \ 0.$$
Then $\alpha + \beta \leq \theta$. 

Let
$$J_+=\left\{j:\ \arg u_j > 0 \right\} \quad \text{and} \quad J_-=\left\{j: \ \arg u_j < 0 \right\}.$$ 
Next, without loss of generality, we assume that $\arg u_j =\alpha$ for all $j \in J_+$ and that $\arg u_j =-\beta$ for all $j \in J_-$.
Indeed, suppose that $\arg u_1 =\alpha_1$ where $0 < \alpha_1< \alpha$. We can modify
$$u_1 \longmapsto {\sin \alpha_1 \over \sin \alpha} e^{i(\alpha-\alpha_1)} u_1$$
(we rotate and shrink $u_1$ so as to make its argument equal to $\alpha$ and leave $\Im\thinspace u_1$ intact).
The sum (2.2.1) gets smaller while all other conditions and the sum (2.2.2) remain intact. Rescaling
$u_j \longmapsto \tau u_j$ for some real $\tau>1$, we make (2.2.1) equal to 1 and increase (2.2.2), while keeping other constraints of the lemma intact.
The case of $\arg u_j > -\beta$ for some $j \in J_-$ is handled similarly. 

Next, without loss of generality, we assume that $\alpha +\beta = \theta$. Indeed, if $\alpha + \beta < \theta$, we can rotate and scale vectors  $u_j$ as above, so that the sum (2.2.2) increases while all other conditions are satisified.

Now, let 
$$u_+=\sum_{j \in J_+} u_j  \quad \text{and} \quad u_-=\sum_{j \in J_-} u_j.$$
Then $\arg u_+=\alpha$, $\arg u_-=-\beta$, $\Im\left(u_+ + u_-\right)=0$, $|u_+| +|u_-|=1$ and (2.2.2) is equal to 
$\left| \Im\thinspace u_+\right| + \left|\Im \thinspace u_-\right|$. 

Denoting $a=|u_+|$ and $b=|u_-|$, we have $a+b=1$ and $a \sin \alpha -b \sin \beta =0$, from which 
$$a={\sin \beta \over \sin \alpha + \sin \beta} \quad \text{and} \quad b={\sin \alpha \over \sin \alpha + \sin \beta}$$
and so 
$$\left| \Im\thinspace u_+ \right| + \left| \Im\thinspace u_-\right| =  {2 \sin \alpha \sin \beta \over \sin \alpha + \sin \beta}.$$
Now, the function 
$$\alpha \longmapsto {1 \over \sin \alpha} \quad \text{for} \quad 0 \leq \alpha \leq {\pi \over 2}$$
is convex and hence the minimum of 
$${\sin \alpha + \sin \beta \over \sin \alpha \sin \beta}={1 \over \sin \alpha} + {1 \over \sin \beta} $$
on the interval $\alpha+\beta=\theta$, $\alpha, \beta \geq 0$, is attained at $\alpha =\beta =\theta/2$.
The proof now follows.
{\hfill \hfill \hfill} \qed
\enddemo

We need another geometric lemma.
\proclaim{(2.3) Lemma} Let $u_1, \ldots, u_n \in {\Bbb C}$ be non-zero complex numbers such that the angle between any two does not exceed $\theta$ for some 
$0 \leq \theta < 2\pi/3$. Let $u=u_1+ \ldots + u_n$. Then 
$$|u| \ \geq \ \left(\cos {\theta \over 2} \right) \sum_{k=1}^n |u_k|.$$
\endproclaim
\demo{Proof} This is Lemma 3.1 of \cite{Ba15} and Lemma 3.6.3 of \cite{Ba16}.
{\hfill \hfill \hfill} \qed
\enddemo 

\head 3. Proof of Theorem 1.3 \endhead

We identify the space of $n \times n$ zero-diagonal complex symmetric matrices $Z=\left(z_{ij}\right)$ with ${\Bbb C}^{n \choose 2}$. Given 
$\delta \geq \eta > 0$, we define a domain $\UU(\delta, \eta) =\UU_{n, m}(\delta, \eta)  \subset {\Bbb C}^{n \choose 2}$ by
$$\UU(\delta, \eta)=\left\{ Z=\left(z_{ij}\right): \ \left| \Re\thinspace z_{ij}\right| \ \leq \ {\delta \over m-1} \quad \text{and} \quad 
\left| \Im\thinspace z_{ij}\right| \ \leq \ {\eta \over m-1} \right\}.$$
We note that the Euclidean distance (in ${\Bbb R}^2={\Bbb C}$) between any two points in $\UU(\delta, \eta)$ does not exceed 
$${\sqrt{(2\delta)^2 + (2 \eta)^2} \over m-1} \ \leq \ {2\sqrt{2} \delta \over m-1}.$$
We will prove by descending induction on $|\Omega|$ that $P_{\Omega}(Z) \ne 0$ for all $Z \in \UU(\delta, \eta)$ and that, moreover, a number of stronger conditions are met. The induction is based on the following two lemmas that describe how $P_{\Omega}(Z)$ changes when only the entries in the $i$-th row and column of $Z$ change. The first lemma deals with the case of $i \in \Omega$.

\proclaim{(3.1) Lemma} Let us fix $\Omega \subset \{1, \ldots, n\}$ such that $ |\Omega| < m$. Suppose that for any 
$Z \in \UU(\delta, \eta)$ and any $j, k \notin \Omega$, we have $P_{\Omega \cup \{j\}}(Z) \ne 0$, 
$P_{\Omega \cup \{k\}}(Z) \ne 0$ and the angle between the two non-zero complex numbers does not exceed $\theta$ for some 
$0 < \theta \leq \pi/2$. 
Then 
\roster
\item We have 
$$P_{\Omega}(Z) \ne 0 \quad \text{for all} \quad Z \in \UU(\delta, \eta).$$
\item Suppose additionally, that $\Omega \ne \emptyset $ and let us fix an $i \in \Omega$. Let $Z', Z'' \in \UU(\delta, \eta)$ be two matrices that differ only in the coordinates
$z_{ij}=z_{ji}$ for $j \ne i$. Then 
$$\left|{P_{\Omega}(Z') \over P_{\Omega}(Z'')}\right| \ \leq \ e^{6\delta}$$ 
and the angle between $P_{\Omega}(Z') \ne 0$ and $P_{\Omega}(Z'') \ne 0$ does not exceed 
$$2 \delta \tan {\theta \over 2} + 5 \eta.$$
\endroster 
\endproclaim
\demo{Proof} It follows from (2.1.1) and Lemma 2.3 that 
$$\left|P_{\Omega}(Z)\right| \ \geq \ {\cos(\theta/2) \over m-|\Omega|} \sum_{j \notin \Omega} \left| P_{\Omega \cup \{j\}}(Z)\right| \ \geq \ {1 \over (m-1) \sqrt{2}} 
 \sum_{j \notin \Omega} \left| P_{\Omega \cup \{j\}}(Z)\right|. 
\tag3.1.1$$
In particular, Part (1) follows. 

To prove Part (2), let us choose a branch of $\ln P_{\Omega}(Z)$ for $Z \in \UU(\delta, \eta)$. 
For $0 \leq t \leq 1$, let $Z(t) = tZ'' +(1-t) Z'$. 
Then 
$$\split \ln P_{\Omega}(Z'') - \ln P_{\Omega}(Z') =& \int_0^1 {d \over dt} \ln  P_{\Omega} \left(Z(t)\right) \ dt\\=
&\int_0^1 \sum_{j: \ j \ne i} \left(z_{ij}'' - z'_{ij}\right) {\partial \over \partial z_{ij}} \ln P_{\Omega}(Z) \Big|_{Z= Z(t)} \ dt.
\endsplit $$
Using (2.1.2), we conclude that 
$${\partial \over \partial z_{ij}} \ln P_{\Omega}(Z) = \cases 1 &\text{if\ } j \in \Omega, \\ P_{\Omega \cup \{j\}}(Z)/P_{\Omega}(Z) 
&\text{if\ } j \notin \Omega, \endcases$$
and hence
$$\aligned \ln P_{\Omega}(Z'') - \ln P_{\Omega}(Z') =& \sum_{j \in \Omega, j \ne i} \left(z_{ij}'' - z_{ij}'\right) \\
&\quad +\int_0^1 \sum_{j \notin \Omega} \left(z_{ij}''- z_{ij}' \right) {P_{\Omega \cup \{j\}}(Z(t)) \over P_{\Omega}(Z(t))}  \ dt. \endaligned \tag3.1.2$$
Using (3.1.1), we get from (3.1.2) that
$$\split \left| \Re\thinspace \ln P_{\Omega}(Z'') - \Re\thinspace \ln P_{\Omega}(Z') \right| \ \leq \ & 2\delta  + {(m-1) \sqrt{2}}\max_{j \ne \Omega} \left| z_{ij}''-z_{ij}'\right| \\ \leq \ &2 \delta + 4 \delta = 6\delta \endsplit $$
and hence 
$$\left| {P_{\Omega}(Z') \over P_{\Omega}(Z'')} \right| \ \leq \ e^{6 \delta},$$
as claimed.

From (2.1.1),  for all $Z \in \UU(\delta, \eta)$ we have that
$$\sum_{j \notin \Omega} {P_{\Omega \cup \{j\}}(Z) \over P_{\Omega}(Z)}= m-|\Omega|$$ 
is real, while from (3.1.1), we conclude that 
$$\sum_{j \notin \Omega} \left| P_{\Omega \cup \{j\}}(Z) \over P_{\Omega}(Z) \right| \ \leq \ {m-|\Omega| \over \cos (\theta/2) } \ \leq \ 
{m-1 \over \cos (\theta/2)}.$$
Applying Lemma 2.2 with $u_j=P_{\Omega \cup \{j\}}(Z)/P_{\Omega}(Z)$, we conclude that 
$$\sum_{j \notin \Omega} \left| \Im\thinspace {P_{\Omega \cup \{j\}}(Z) \over P_{\Omega}(Z)} \right| \ \leq \ (m-1) \tan {\theta \over 2}.$$
Therefore, from (3.1.2),
$$\split \left| \Im\thinspace \ln P_{\Omega}(Z'') - \Im \ln P_{\Omega}(Z') \right| \ \leq \ &2 \eta + 
(m-1) \tan {\theta \over 2} \max_{j \notin \Omega} \left|\Re\thinspace z_{ij}'' - \Re\thinspace z_{ij}' \right| \\ &\quad +(m-1)\sqrt{2} \max_{j \notin \Omega} \left|\Im\thinspace z_{ij}'' - \Im\thinspace z_{ij}' \right| \\ &\leq \ 2 \delta \tan {\theta \over 2} + 5 \eta. \endsplit$$
Hence the angle between $P_{\Omega}(Z'')$ and $P_{\Omega}(Z')$ does not exceed $2 \delta \tan {\theta \over 2} + 5\eta$, as claimed.
{\hfill \hfill \hfill} \qed
\enddemo

The second lemma shows that  $P_{\Omega}(Z)$ does not change much if only the entries of $Z$ in the $i$-th row and column are changed for some $i \notin \Omega$,
assuming that $n \gg m$.

\proclaim{(3.2) Lemma} Let us fix an $\Omega \subset \{1, \ldots, n\}$, $|\Omega| \leq m-1$. 
Suppose for any 
$i,j \notin \Omega$ and all $Z \in \UU(\delta, \eta)$ we have $P_{\Omega \cup \{i\}}(Z) \ne 0$, $P_{\Omega \cup \{j\}}(Z) \ne 0$ and the angle between the 
two complex numbers does not exceed $\pi/2$ and that 
$$\left| {P_{\Omega \cup \{i\}}(Z) \over P_{\Omega \cup \{j\}}(Z)} \right| \ \leq \ \lambda$$
for some $\lambda \geq 1$.

In addition, suppose that if $|\Omega| \leq m-2$ then for any distinct $i, j, k \notin \Omega$ and all 
$Z \in \UU(\delta, \eta)$ we have $P_{\Omega \cup \{i, j\}}(Z) \ne 0$, $P_{\Omega \cup \{i, k\}}(Z) \ne 0$ and the angle between the two complex numbers 
does not exceed $\pi/2$.

Let us fix an $i \notin \Omega$ and let $Z', Z'' \in \UU(\delta, \eta)$ be two matrices that differ only in the coordinates $z_{ij}=z_{ji}$ for 
$j \ne i$. Then 
$$\left| {P_{\Omega}(Z') \over P_{\Omega}(Z'')} \right| \ \leq \ \exp\left\{{10 \delta \lambda m \over n-1} \right\}$$
and the angle between $P_{\Omega}(Z') \ne 0$ and $P_{\Omega}(Z'') \ne 0$ does not exceed 
$${10\delta \lambda m \over n-1}.$$
\endproclaim
\demo{Proof} It follows from Lemma 3.1 that $P_{\Omega}(Z) \ne 0$ for all $Z \in \UU(\delta, \eta)$.

 Arguing as in the proof of Lemma 3.1, we introduce $Z(t)=t Z''+(1-t) Z'$ and write 
$$\ln P_{\Omega}(Z'')-\ln P_{\Omega}(Z')=\int_0^1 \sum_{j:\ j \ne i} \left(z_{ij}''-z_{ij}'\right) {\partial \over \partial z_{ij}} \ln P_{\Omega}(Z)\Big|_{Z=Z(t)} \ dt.$$
From (2.1.2), we write 
$$\aligned \ln P_{\Omega}(Z'')-\ln P_{\Omega}(Z')=&\int_0^1 \sum_{j \in \Omega} \left(z_{ij}''-z_{ij}'\right) {P_{\Omega \cup \{i\}}\bigl(Z(t)\bigr) \over P_{\Omega}\bigl(Z(t)\bigr)} \\ &\quad +
\sum_{j \notin \Omega, j \ne i}  \left(z_{ij}''-z_{ij}'\right) {P_{\Omega \cup \{i, j\}}\bigl(Z(t)\bigr) \over P_{\Omega}\bigl(Z(t)\bigr)} \ dt. \endaligned \tag3.2.1$$

Suppose first that $|\Omega| \leq m-2$. From (2.1.1), we have 
$$P_{\Omega \cup \{i\}}(Z)={1 \over m-|\Omega| -1} \sum_{j \notin \Omega, j \ne i} P_{\Omega \cup \{i, j\}}(Z).$$
Applying Lemma 2.3, we get that
$$\sum_{j \notin \Omega, j \ne i} \left| P_{\Omega \cup \{i, j\}}(Z) \right| \ \leq \ (m-1)\sqrt{2} \left| P_{\Omega \cup \{i\}}(Z) \right| \tag3.2.2$$
for all $Z \in \UU(\delta, \eta)$.

Since by (2.1.1) we also have 
$$P_{\Omega}(Z) ={1 \over m-|\Omega|} \sum_{j \notin \Omega} P_{\Omega \cup \{j\}}(Z),$$
applying Lemma 2.3, we conclude that 
$$\sum_{j \notin \Omega} \left| P_{\Omega \cup \{j\}}(Z)\right| \ \leq \ (m-|\Omega|) \sqrt{2} \left| P_{\Omega}(Z)\right|.$$
Hence for all $i \notin \Omega$, we have
$$\left|P_{\Omega \cup \{i\}}(Z)\right| \ \leq \ {\lambda(m-|\Omega|) \sqrt{2} \over n-|\Omega|} \left|P_{\Omega}(Z)\right| \ \leq \ {\lambda m \sqrt{2} \over n} 
\left|P_{\Omega}(Z)\right|. \tag3.2.3$$
Combining (3.2.3) and (3.2.2), we get 
$$\sum_{j \notin \Omega, j \ne i} \left| P_{\Omega \cup \{i, j\}}(Z)\right| \ \leq \ {2\lambda m(m-1) \over n} \left|P_{\Omega}(Z)\right|. \tag3.2.4$$
Combining (3.2.1), (3.2.2), (3.2.3) and (3.2.4), we get 
$$\split \left| \ln P_{\Omega}(Z'') - \ln P_{\Omega}(Z') \right| \ \leq \ & {2 \sqrt{2} \delta \over m-1} \cdot {\lambda |\Omega| (m-|\Omega|) \sqrt{2} \over n-|\Omega|}  +
{2 \sqrt{2} \delta \over m-1} \cdot {2 \lambda m (m-1) \over n} \\ \leq \ & {4 \delta \lambda m \over n-1} + {4 \sqrt{2} \delta \lambda m \over n} \ \leq \ {10 \delta \lambda m \over n-1} . \endsplit$$
If $|\Omega| =m-1$ then from (3.2.1) and (3.2.3), we get 
$$\left| \ln P_{\Omega}(Z'') -\ln P_{\Omega}(Z')\right| \ \leq \ {2 \sqrt{2} \delta \over m-1} \cdot {\lambda m \sqrt{2} \over n} \ \leq \ {4 \delta \lambda m \over n-1},$$
which concludes the proof.
{\hfill \hfill \hfill} \qed
\enddemo

Now we are ready to prove Theorem 1.3.

\subhead (3.3) Proof of Theorem 1.3 \endsubhead Given $0 < \delta < 1$, we choose 
$0 < \theta < \pi/2$ so that 
$$2 \delta \tan {\theta \over 2} \ < \ \theta.$$
We then choose $\eta > 0$ such that 
$$2 \delta \tan {\theta \over 2} + 5\eta \ < \ \theta.$$
We choose 
$$\lambda > e^{6 \delta}$$ and 
choose $\omega >1$ so that 
$$2 \delta \tan {\theta \over 2} + 5\eta + {10 \delta \lambda m  \over n-1}  \ \leq \ \theta \quad \text{and} \quad
 \exp\left\{ 6\delta+ {10 \delta \lambda m \over n-1}\right\} \leq \lambda$$
 whenever $n \geq \omega m$. 

Suppose that $n \geq \omega m$.
We prove by descending induction on $r=m, m-1, \ldots, 1$ that if $\Omega_1, \Omega_2 \in \{1, \ldots, n\}$ are two sets such that 
$|\Omega_1|=|\Omega_2|=r$ and $\left|\Omega_1 \Delta \Omega_2\right|=2$ then for all $Z \in \UU(\delta, \eta)$ we have $P_{\Omega_1}(Z) \ne 0$, $P_{\Omega_2}(Z) \ne 0$, the 
angle between $P_{\Omega_1}(Z)$ and $P_{\Omega_2}(Z)$ does not exceed $\theta$ while the ratio of $|P_{\Omega_1}(Z)|$ and $|P_{\Omega_2}(Z)|$
does not exceed $\lambda$.

Assume that $r=m$. Without loss of generality, we assume that $\Omega_1=\Omega \cup \{1\}$ and $\Omega_2=\Omega \cup \{2\}$ for some 
$\Omega \subset \{3, \ldots, n\}$ such that $|\Omega|=m-1$.
We have 
$$\split P_{\Omega_1}(Z) =&\exp\left\{\sum_{\{i, j\} \subset \Omega} z_{ij}\right\} \exp\left\{\sum_{i \in \Omega} z_{1i}\right\} \quad \text{and} \\
P_{\Omega_2}(Z) =&\exp\left\{\sum_{\{i, j\} \subset \Omega} z_{ij}\right\} \exp\left\{\sum_{i \in \Omega} z_{2i}\right\}. \endsplit$$
Clearly, $P_{\Omega_1}(Z) \ne 0$, $P_{\Omega_2}(Z) \ne 0$, the angle between $P_{\Omega_1}(Z)$ and $P_{\Omega_2}(Z)$ does not exceed
$2 \eta \leq \theta$ while the ratio of $|P_{\Omega_1}(Z)|$ and $|P_{\Omega_2}(Z)|$ does not exceed $e^{2 \delta} \leq \lambda$.

Suppose now that the statements hold for all subsets $\Omega \subset \{1, \ldots, n\}$ of cardinality at least $r+1$ for some $r \leq m-1$ and let $\Omega_1, \Omega_2 \subset \{1, \ldots, n\}$
we two subsets of cardinality $r \geq 1$ such that $|\Omega_1 \Delta \Omega_2| =2$. Again, without loss of generality, we assume that 
$\Omega_1 =\Omega \cup \{1\}$ and $\Omega_2=\Omega \cup \{2\}$ for some $\Omega \subset \{3, \ldots, n\}$ such that $|\Omega|=r-1$.
Then we observe that $P_{\Omega_2}(Z)=P_{\Omega_1}(Z')$, where 
$$z'_{1i}=z'_{i1}=z_{2i}=z_{i2} \quad \text{and} \quad z_{2i}'=z_{i2}'=z_{1i}=z_{i1} \quad \text{for} \quad i \ne 1, 2,$$
while all other entries of $Z$ and $Z'$ coincide. Applying Lemma 3.1 and Lemma 3.2 and the induction hypothesis to sets $\Omega_1 \cup \{j\}$ for 
$j \notin \Omega_1$ and $\Omega_1 \cup\{j, k\}$ for $j, k\notin \Omega_1$,
 we conclude that the angle between $P_{\Omega_1}(Z) \ne 0$ and $P_{\Omega_2}(Z) \ne 0$ does not exceed 
$$2 \delta \tan{\theta \over 2} + 5 \eta + {10 \delta \lambda m \over n-1} \ \leq \ \theta,$$
while the ratio of $|P_{\Omega_1}(Z)|$ and $|P_{\Omega_2}(Z)|$ does not exceed
$$ \exp\left\{ 6 \delta+{10 \delta \lambda m \over n-1}\right\} \leq \lambda.$$

This proves that $P_{\{i\}}(Z) \ne 0$ for all $i \in \{1, \ldots, n\}$ and all $Z \in \UU(\delta, \eta)$ and that the angle between 
$P_{\{i\}}(Z) \ne 0$ and $P_{\{j\}}(Z) \ne 0$ does not exceed $\theta$ for all $i, j\in \{1, \ldots, n\}$. From (2.1.1) we conclude that 
$P_m(Z) =P_{\emptyset}(Z) \ne 0$ for all $Z \in \UU(\delta, \eta)$.
{\hfill \hfill \hfill} \qed

\head 4. Computing the partition function \endhead

Here we show how to compute the density partition function $\den_m(G; \gamma)$. First, we make a change of coordinates 
to convert the partition function $P_m(Z)$ of Section 1.2 into a multivariate polynomial.
\subhead (4.1) A polynomial version of $P_m(Z)$ \endsubhead For an $n \times n$ complex symmetric matrix $W=\left(w_{ij}\right)$ with zero diagonal, we 
define 
$$p_m(W)={n \choose m}^{-1}\sum\Sb S \subset \{1, \ldots, n\} \\ |S|=m \endSb \prod\Sb \{i, j\} \subset S \\ i \ne j \endSb \left(1+w_{ij}\right).$$
Hence $p_m(W)$ is a polynomial of degree ${m \choose 2}$ in the entries $w_{ij}$ and, assuming that $|w_{ij}| < 1$ for all $i, j$, we can write 
$$p_m(W)={n \choose m} P_m(Z) \quad \text{where} \quad Z=\left(z_{ij} \right) \quad \text{and} \quad z_{ij} =\ln \left(1+ w_{ij}\right)$$
(we choose the standard branch of the logarithm in the right half-plane of ${\Bbb C}$).
Theorem 1.3 implies that for every $0 < \delta < 1$ there is $\eta=\eta(\delta) >0$ and $\omega=\omega(\delta) >1$ such that 
$$\aligned p_m(W) \ne 0 \quad \text{whenever} \quad &\left| \Re\thinspace \ln \left(1+w_{ij}\right)\right| \ \leq \ {\delta \over m-1}, \\ &\left| \Im\thinspace \ln \left(1+w_{ij}\right)\right| 
\ \leq \ {\eta \over m-1} \quad \text{and}\\  &n \geq \omega m. \endaligned \tag4.1.1$$
To compute $\den_m(G; \gamma)$ for a given $0 < \gamma <1$ and a given graph $G=(V, E)$, we define
$$w_{ij}=\cases \exp\left\{{\gamma \over m-1}\right\}-1 &\text{if\ } \{i, j\} \in E, \\ \exp\left\{-{\gamma \over m-1}\right\}-1 &\text{if\ } \{i, j\} \notin E. \endcases \tag4.1.2$$
Then, by (1.2.2), we have 
$$\den_m(G; \gamma)= \exp\left\{ {\gamma m \over 2}\right\} p_m(W). \tag4.1.3$$

The interpolation method is based on the following simple lemma.
\proclaim{(4.2) Lemma} Let $g: {\Bbb C} \longrightarrow {\Bbb C}$ be a univariate polynomial and suppose that $g(z) \ne 0$ provided $|z| < \beta$
where $\beta > 1$ is some real number. Let us choose a branch of $f(z)=\ln g(z)$ in the disc $|z| < \beta$ and let 
$$T_r(z)=f(0)+ \sum_{k=1}^r {f^{(k)}(0) \over k!} z^k$$
be the Taylor polynomial of $f$ of degree $r$ computed at $z=0$. Then 
$$\left| f(1) - T_r(1)\right| \ \leq \ {\deg g \over \beta^r(\beta-1) (r+1)}.$$
\endproclaim
\demo{Proof} This is Lemma 2.2.1 of \cite{Ba16}, see also Lemma 1.1 of \cite{Ba15}.
{\hfill \hfill \hfill} \qed
\enddemo
The gist of Lemma 4.2 is that to approximate $f(1)$ within an additive error $\epsilon$, it suffices to compute the Taylor polynomial of $f(z)$ at $0$ of degree 
$r =O_{\beta} \left( \ln \deg g - \ln \epsilon \right)$, where the implicit constant in the ``$O$" notation depends on $\beta$ alone.
We would like to apply Lemma 4.2 to the univariate polynomial 
$$h(z)={n \choose m}^{-1} \sum\Sb S \subset \{1, \ldots, n\} \\ |S|=m \endSb \prod\Sb \{i, j\} \subset S \\ i \ne j \endSb \left( 1+ z w_{ij}\right), \tag4.2.1$$
where $w_{ij}$ are defined by (4.1.2). Indeed, the value we are ultimately interested is $h(1)=p_m(W)$. However, Lemma 4.2 requires that $h(z) \ne 0$ in a disc of some radius 
$\beta >1$, whereas (4.1.1) only guarantees that $h(z) \ne 0$ for $z$ in a neighborhood of the interval $[0, 1] \subset {\Bbb C}$. To remedy this, we compose 
$h$ with a polynomial $\phi: {\Bbb C} \longrightarrow {\Bbb C}$ such that $\phi(0)=0$, $\phi(1)=1$ and $\phi$ maps the disc $|z| < \beta$ for some $\beta >1$ inside the prescribed neighborhood of $[0, 1] \subset {\Bbb C}$. We then apply Lemma 4.2 to the composition $g(z)=h((\phi(z))$. The following lemma provides an explicit construction of $\phi$.

\proclaim{(4.3) Lemma} For $0 < \rho < 1$, we define 
$$\split &\alpha=\alpha(\rho)=1-e^{-{1 \over \rho}}, \quad \beta=\beta(\rho)={1 - e^{-1-{1 \over \rho}} \over 1-e^{-{1 \over \rho}}} \ > \ 1, \\
&N=N(\rho)=\left\lfloor \left(1+ {1 \over \rho}\right) e^{1+ {1 \over \rho}}\right\rfloor, \quad \sigma=\sigma(\rho)=\sum_{k=1}^N {\alpha^k \over k} \quad \text{and} \\
&\phi(z)=\phi_{\rho}(z)={1 \over \sigma} \sum_{k=1}^N {(\alpha z)^k \over k}. \endsplit$$
Then $\phi: {\Bbb C} \longrightarrow {\Bbb C}$ is a polynomial of degree $N$ such that $\phi(0)=0$, $\phi(1)=1$,
$$-\rho \ \leq \ \Re\thinspace \phi(z) \ \leq \ 1+ 2\rho \quad \text{and} \quad \left| \Im\thinspace \phi(z)\right| \ \leq \ 2 \rho$$
provided $|z| \leq \beta$.
\endproclaim 
\demo{Proof} This is Lemma 2.2.3 of \cite{Ba16}.
{\hfill \hfill \hfill} \qed
\enddemo

Lemma 4.2 also requires the derivatives $f^{(k)}(0)$ of $f(z)=\ln g(z)$ at $z=0$. Those, however, can be easily computed from the derivatives $g^{(k)}(0)$, as described in Section 2.2.2 of \cite{Ba16}, see also Section 2.1 of \cite{Ba15}. We briefly sketch how. 

\subhead (4.4) Computing derivatives \endsubhead Suppose that $f(z)=\ln g(z)$ as in Lemma 4.2. Then 
$$f'(z)={g'(z) \over g(z)} \quad \text{and} \quad g'(z)=f'(z) g(z).$$
Differentiating the product $k-1$ times, we obtain 
$$g^{(k)}(0)=\sum_{j=0}^{k-1} {k-1 \choose j} f^{(k-j)}(0) g^{(j)}(0) \quad \text{for} \quad k=1, \ldots, r. \tag4.4.1$$
We interpret (4.4.1) as a system of linear equations in variables $f^{(k)}(0)$ for $k=1, \ldots, r$ with coefficients $g^{(k)}(0)$ for $k=0, \ldots, r$. This is a triangular system of linear equations with non-zero entries $g^{(0)}(0)=g(0)$ on the diagonal, that can be solved in $O(r^2)$ time, provided the values of $g^{(k)}(0)$ are known.

To supply the last ingredient of the algorithm, we show how to compute $h^{(k)}(0)$ for $k=0, \ldots, r$, where $h$ is the polynomial defined by (4.2.1). This is also done in \cite{Ba15}, but we reproduce it here for completeness. 

We have 
$$h^{(k)}(0)={n \choose m}^{-1} \sum\Sb S \subset \{1, \ldots, n\} \\ |S|=m \endSb \sum_{\{i_1, j_1\}, \ldots, \{i_k, j_k\} \subset S} 
w_{i_1 j_1} \cdots  w_{i_k j_k},$$
where the inner sum is taken over all ordered collections of distinct unordered pairs $\{i_1, j_1\}, \ldots, \{i_k, j_k\} \subset S$.
For such a collection, say $I$, let $\nu(I)$ be the number of distinct vertices among $i_1, j_1, \ldots, i_k, j_k$. Then there are exactly ${n-\nu(I) \choose m-\nu(I)}$ different $m$-subsets $S$ containing the edges from $I$ and we can rewrite the above sum
as 
$$h^{(k)}(0) ={n \choose m}^{-1} \sum_{I=\left( \{i_1, j_1\}, \ldots, \{i_k, j_k\}\right)} {n - \nu(I) \choose m-\nu(I)} w_{i_1 j_1} \cdots w_{i_k j_k}, \tag4.4.2$$
where the sum is taken over all ordered collections of $k$ distinct edges in $G$. It is clear now that $h^{(k)}(0)$ can be computed in \
$n^{O(k)}$ time by the exhaustive enumeration of all possible collections of $k$ edges. 

In Section 5 we present faster formulas for computing $h^{(2)}(0)$ and $h^{(3)}(0)$ that we used for our numerical experiments.

\subhead (4.5) The algorithm \endsubhead Let us fix $0 < \gamma < 1$. Below we summarize the algorithm for computing 
$\den_m(G; \gamma)$ within relative error $0< \epsilon < 1$, by which we understand computing $\ln \den_m(G; \gamma)$ within 
additive error $\epsilon$. We assume that $m \geq 4$ and that 
$n \geq \omega m$ for some $\omega=\omega(\gamma) > 1$, to be specified below. 

Given a graph $G=(V, E)$ with set $V=\{1, \ldots, n\}$ of vertices, and an integer $m \leq n$,  we compute the $n \times n$ symmetric matrix $W=\left(w_{ij}\right)$ by (4.1.2). 
Since $m \geq 4$, we have $|w_{ij}| \leq 0.4$ for all $i, j$.

Our goal is to compute $p_m(W)=h(1)$, where $h$ is the univariate polynomial defined by (4.2.1). We note that $\deg h = {m \choose 2}$.
 
 Let us choose $1 > \delta > \gamma$ and let $\eta=\eta(\delta) >0$ and $\omega=\omega(\delta) >1$ be the numbers of Theorem 1.3 and in (4.1.1).
 We find $\rho=\rho(\delta) >0$ such that 
 $$\left| \Re\thinspace \ln \left(1 + z w_{ij}\right) \right| \ \leq \ {\delta \over m-1} \quad \text{and} \quad 
 \left| \Im\thinspace \ln \left(1+ zw_{ij}\right) \right| \ \leq \ {\eta \over m-1}$$
 as long as 
 $$-\rho \ \leq \ \Re\thinspace z \ \leq \ 1+ \rho \quad \text{and} \quad  \left| \Im\thinspace z \right| \ \leq \ \rho. \tag4.5.1$$
 Indeed, if $z \in [0, 1]$ then 
 $$-{\gamma \over m-1} \ \leq \ \ln \left(1 +z w_{ij}\right) \ \leq \ {\gamma \over m-1}$$
 and for $|z| \leq 2$, we have
 $$\left|{d \over dz} \ln \left(1+ z w_{ij}\right) \right| =\left| {w_{ij} \over 1+z w_{ij}} \right| \ \leq \ {10 \over m-1}$$
so the desired $\rho$ can indeed be found.
 
 It follows by (4.1.1) that $h(z) \ne 0$ as long as $n \geq \omega m$ and (4.5.1) holds. 
  
Using Lemma 4.3, we construct a polynomial $\phi: {\Bbb C} \longrightarrow {\Bbb C}$ of some degree $N=N(\rho)=N(\delta)$ such that 
$\phi(0)=0$, $\phi(1)=1$ and 
$$-\rho \ \leq \ \Re\thinspace \phi(z) \ \leq \ 1+ \rho \quad \text{and} \quad  \left| \Im\thinspace \phi(z) \right| \ \leq \ \rho$$
 as long as $|z| \leq \beta$ for some $\beta=\beta(\rho) =\beta(\delta)>1$. 
We define 
$$g(z) = h(\phi(z))$$ 
and our goal is to compute $g(1)=h(\phi(1))$. We note that 
$$\deg g \ \leq \ N \deg h = N {m \choose 2}.$$ We choose a branch of $f(z)=\ln g(z)$ for $z$ satisfying (4.5.1).

 Using Lemma 4.2, we find an integer $r=O_{\rho}\left( \ln m - \ln \epsilon\right)=O_{\delta}\left(\ln m - \ln \epsilon\right)$ 
such that $$|T_r(1) - f(1)| \ \leq \ \epsilon,$$ where $T_r(z)$ is the Taylor polynomial of $f(z)$ of degree $r$, computed at $z=0$. 
The implicit constant in the ``$O$" notation depends only on $\rho$, which in turn depends only on $\delta$. Hence our goal is to compute $T_r(1)$, for which we need to compute $f^{(k)}(0)$ for $k=1, \ldots, r$. As in Section 4.4, we reduce it in $O(r^2)$ time to computing 
$g^{(k)}(0)$ for $k=1, \ldots, r$. Note that 
$$g(0)=h(\phi(0))= h(0) =1.$$
Let $\phi_r(z)$ be the truncation of the polynomial $\phi(z)$ obtained by discarding all monomials of degree higher than $r$. 
Similarly, let $h_r(z)$ be the truncation of the polynomial $h(z)$, obtained by discarding all monomial of degree higher than $r$. 
We compute $h_r(z)$ as in Section 4.4 in $n^{O(r)}$ time. Finally, we compute the truncation of the composition $h_r(\phi_r(z))$. A fast (polynomial in $r$) way to do it, is to use Horner's method: assuming that 
$$h_r(z)=\sum_{k=0}^r b_k z^k,$$ 
we successively compute
$$\split &b_r \phi_r(z) + b_{r-1}, \quad \left( b_r \phi_r(z) + b_{r-1}\right) \phi_r(z)+ b_{r-2}, \\
&\left( \left( b_r \phi_r(z) + b_{r-1}\right) \phi_r(z)+ b_{r-2}\right) \phi_r(z) + b_{r-3}, \ldots \endsplit$$
discarding on the way all monomials of degree higher than $r$. In the end, we have computed $g^{(k)}(0)$ for $k=0, \ldots, r$ and 
hence $f^{(k)}(0)$ for $k=0, \ldots, r$ and hence $T_m(1)$ approximating $f(1)=\ln h(1)$ within additive error $\epsilon$. From (4.1.3), we compute 
$$\den_m(G; \gamma) = \exp\left\{ { \gamma m \over 2} \right\} h(1)$$ 
within relative error $\epsilon > 0$.

\head 5. Remarks on the practical implementation \endhead

We implemented a {\it much} simplified version of the algorithm. Given a graph $G=(V, E)$ with set $V=\{1, \ldots, n\}$ of vertices and an integer $2 \leq m \leq n$, we define 
the $n \times n$ matrix $=\left(w_{ij}\right)$ by 
$$w_{ij}=\cases \alpha &\text{if\ } \{i, j\} \in E \\ -\alpha &\text{if\ } \{i, j\} \notin E, \endcases$$
where $0 < \alpha < 1$ is a parameter.

We the consider the polynomial $h(z)$ defined by (4.2.1) and let $f(z)=\ln h(z)$. 

Our goal is to approximate $f(1)$ and hence 
$$\split h(1) = &\sum\Sb S \subset \{1, \ldots, n \} \\ |S| =m \endSb (1+ \alpha)^{{m \choose 2} \sigma(S)} (1- \alpha)^{{m \choose 2} (1- \sigma(S))} \\=
&(1-\alpha)^{{m \choose 2}} \den_m(G; \gamma), \quad \text{where} \quad \gamma={m -1 \over 2} \ln {1+\alpha \over 1-\alpha}. \endsplit$$  
We approximate $f(1)$ by the degree $r$ Taylor polynomial of $f(z)$ computed at $z=0$. The results of \cite{Ba15} suggest that for $\alpha=O\left(1/m\right)$, we should get a reasonable approximation if we use $r \sim \ln m$. The results of our numerical experiments suggest that we get reasonable approximations if we use $\alpha=\Omega(1)$ and 
$r=2$ or $r=3$. In short, on the examples we tested, the quality of approximation was more consistent with the quality of the Taylor polynomial approximation of $\ln (1\pm \alpha)$.

We provide below the explicit formulas for the approximations up to degree 3, in case the reader will be interested to do some numerical experiments. We interpret $w_{ij}$ as weights on the edges of a complete graph with $n$ vertices. Borrowing an idea from \cite{PR17}, we express the derivatives $f^{(k)}(0)$ in terms of various sums associated with {\it connected} subgraphs, since it improves the computational complexity of the algorithm.

It is convenient to introduce the following sums:
$$A_1=\sum_{\{i, j\}} w_{ij},$$
where the sum is taken over all unordered pairs $\{i, j\}$ of distinct indices;
 $$B_1=\sum_{\{i, j\}} w_{ij}^2, \quad B_2=\sum_{j, \{i, k\}} w_{ij}w_{jk},$$
 where in the formula for $B_1$ the sum is taken oven all unordered pairs $\{i, j\}$ of distinct indices and in $B_2$ the sum is taken over all pairs consisting of an index $j$ and 
 an unordered pair $\{i, k\}$, so that all three indices are distinct; and 
 $$\split &C_1 = \sum_{\{i, j\}} w_{ij}^3, \quad C_2 =\sum_{(i,j,k)} w_{ij}^2 w_{jk}, \quad C_3=\sum_{\{i,j,k\}} w_{ij} w_{jk} w_{ki}, \\
 &C_4=\sum_{(i,j,k,l)} w_{ij}w_{jk}w_{kl}, \quad C_5= \sum_{\{j,k,l\}, i} w_{il}w_{ij}w_{ik}, \endsplit$$
where in $C_1$ the sum is taken over all unordered pairs $\{i, j\}$ of distinct indices, in $C_2$ the sum is taken over all ordered triples $(i, j, k)$ of distinct indices, in $C_3$ the sum is taken over all unordered triples of distinct integers, in 
$C_4$, the sum is taken over all ordered $4$-tuples $(i, j, k, l)$ of distinct indices, and in $C_5$ the sum is taken over all pairs consisting of an index $i$ and an unordered triple 
$\{j, k, l\}$ so that all four indices $\{i, j, k, l\}$ are distinct.

\subhead (5.1) First-order approximation \endsubhead Clearly, $h(0)=1$. From (4.4.2), we have 
$$h'(0)={n \choose m}^{-1} {n-2 \choose m-2} \sum_{\{i, j\} \subset \{1, \ldots, n\}} w_{ij}
={m(m-1) \over n(n-1)} A_1.$$
Since $f(0)=\ln h(0)=0$ and $f'(0)=h'(0)/h(0)=h'(0)$, we obtain the first order approximation
$$f(1) \approx h_1(0),$$
where $h_1(0)$ is defined as above. The complexity of computing the first order approximation in $O(n^2)$.

\subhead (5.2) Second-order approximation \endsubhead From (4.4.2), we have 
$$h''(0)={n \choose m}^{-1} \sum_{I=(\{i_1, j_1\}, \{i_2, j_2\})} {n-\nu(I) \choose m-\nu(I)} w_{i_1 j_1} w_{i_2 j_2}.$$
Here $\nu(I)=4$ if the pairs $\{i_1, j_1\}$ and $\{i_2, j_2\}$ are pairwise disjoint and $\nu(I)=3$ if they share exactly one index.
Hence we can write
$$\split h''(0)=& {n \choose m}^{-1} \left(2 {n-3 \choose m-3}  B_2 + {n-4 \choose m-4}\left( A_1^2- 2B_2 - B_1\right)\right)\\
=&2{m(m-1)(m-2) \over n(n-1)(n-2)} B_2 + {m (m-1)(m-2)(m-3) \over n(n-1)(n-2)(n-3)} \left(A_1^2 - 2B_2 -B_1\right). \endsplit$$
Since 
$$f''(0)=h''(0)-\left(h'(0)\right)^2,$$
we obtain the second order approximation:
$$f(1) \approx f'(0) + {1 \over 2} f''(0)= h'(0)-{1 \over 2}\left(h'(0)\right)^2 + {1 \over 2} h''(0),$$
where $h'(0)$ and $h''(0)$ are defined as above. The complexity of computing the second order approximation is $O(n^3)$.

\subhead (5.3) Third-order approximation \endsubhead From (4.4.2), one can deduce that 
$$\split  h'''(0)= &6 {m(m-1)(m-2) \over  n(n-1)(n-2)} C_3 +{m(m-1)(m-2)(m-3) \over n(n-1)(n-2)(n-3)} \left( 6 C_5 + 3C_4\right)  \\&\quad + 6 {m(m-1)(m-2)(m-3)(m-4) \over n(n-1)(n-2)(n-3)(n-4)}
\left( A_1 B_2 -3C_5 -3C_3-C_4-C_2\right) \\
&\qquad +{m(m-1)(m-2)(m-3)(m-4)(m-5) \over n(n-1)(n-2)(n-3)(n-4)(n-5)} \Bigl(A_1^3 +12C_3  -6A_1 B_2 \\ &\qquad \qquad  \qquad \qquad +12 C_5 +3C_4 +6C_2 - 3A_1 B_1+2C_1\Bigr). \endsplit$$
 Since we have 
 $$f'''(0)=h'''(0)-2f''(0)h'(0)-f'(0)h''(0)=2(h'(0))^3-3h'(0)h''(0)+h'''(0),$$
 we obtain the third order approximation approximation
 $$\split f(1) \approx &f'(0) + {1 \over 2} f''(0) + {1 \over 6} f'''(0)\\=&
 h'(0)-{1 \over 2} (h'(0))^2 +{1 \over 2} h''(0) + {1 \over 3} (h'(0))^3 - {1 \over 2} h'(0) h''(0) + {1 \over 6} h'''(0). \endsplit$$
 The complexity of computing the third order approximation is $O\left(n^4\right)$.

\head 6. Proof of Theorem 1.4 \endhead 

We got the idea of the proof from \cite{EM17}, where a similar question about complex zeros of the permanents of matrices with independent random entries was treated.

Applying Jensen's formula, see for example, Section 5.3 of \cite{Ah78},
we obtain
$$\ln \left| h_W(0)\right| =\sum_{s=1}^N \ln {\left| a_{s,W}\right| \over r}  +{1 \over 2\pi} \int_0^{2\pi} \ln \left|  h_W\left( re^{ i\theta}\right)\right| \ d \theta, \tag6.1$$
where $a_{s,W}, s=1, \ldots, N$ are the roots of the polynomial $h_W(z)$ in the disc $|z| < r$ and we assume that $h_W(z)$ has no zeros on the circle $|z|=r$ (since there are only finitely many values of $r$ with roots on the circle $|z|=r$, this assumption is not restrictive).
We have 
$$\ln \left| h_W(0)\right| =0$$
and furthermore, applying Jensen's inequality twice, we bound: 
$$\aligned {1 \over 2\pi} \int_0^{2\pi} \ln \left|  h_W\left( re^{i \theta}\right)\right| \ d \theta \ \leq \ &\ln \left({1 \over 2\pi}  \int_0^{2\pi} \left|  h_W\left( re^{i \theta}\right)\right| \ d \theta  \right) \\
\ \leq \ &{1 \over 2} \ln \left({1 \over 2\pi} \int_0^{2\pi} \left|  h_W\left( re^{i \theta}\right)\right|^2 \ d \theta  \right). \endaligned \tag6.2$$
For a fixed $\theta \in [0, 2\pi]$, we compute the expectation 
$$\split \EE \left|h_W\left(r e^{i \theta}\right)\right|^2=&{n \choose m}^{-2} \sum\Sb S_1, S_2 \subset \{1, \ldots, n\}\\ |S_1|=|S_2|=m \endSb \EE \Biggl(
\prod_{\{j, k\} \subset S_1} \left(1+ re^{i \theta} w_{jk}\right) \\ &\qquad \qquad \qquad  \times  \prod_{\{j, k\} \subset S_2} \left( 1+ re^{-i \theta} w_{jk} \right) \Biggr) \\
=& {n \choose m}^{-2} \sum\Sb S_1, S_2 \subset \{1, \ldots, n\} \\ |S_1| =|S_2|=m \endSb \left(1+r^2\right)^{|S_1 \cap S_2| \choose  2}. \endsplit $$
A subset $S \subset \{1, \ldots, n\}$ of cardinality $l=|S| \leq m$ can be represented as the intersection $S=S_1 \cap S_2$ of $m$-subsets $S_1, S_2$ in 
${n-l \choose m-l} {n-m \choose m-l}$ ways. 
Hence 
$$\EE \left|h_W\left(r e^{i \theta}\right)\right|^2 ={n \choose m}^{-2}\sum_{l=0}^m {n \choose l}{n-l \choose m-l} {n-m \choose m-l} \left(1+r^2\right)^{l \choose 2}. \tag6.3$$
To bound (6.3), we consider the ratio of the $(l+1)$-st term to the $l$-th term:
$$\split &{n-l \over l+1} \cdot {m-l \over n-l} \cdot {m-l \over n-2m+l+1} \cdot \left(1+r^2\right)^l= {(m-l)^2 \left(1+ r^2 \right)^l \over (l+1)(n-2m+l+1)} \\
&\qquad \leq \ {m^2 (1+r^2)^m \over n-2m+1}. \endsplit$$
In particular, if 
$$n \ \geq \ 2m^2 (1+r^2)^m +2m, \tag6.4$$
the ratio does not exceed $1/2$ and hence we can bound the sum (6.3) by 
$$\EE \left| h_W \left(r e^{i\theta}\right) \right|^2 \ \leq \ 2{n \choose m}^{-2} {n \choose m} {n-m \choose m} \ \leq \ 2.$$
Integrating over $\theta$, we conclude that if (6.4) holds then 
$$\EE \left( {1 \over 2\pi}\int_0^{2\pi} \left| h_W\left(r e^{i \theta}\right) \right| \ d \theta \right) \ \leq \ 2.$$
By the Markov inequality, for any $\tau \geq 1$, we get 
$$\Pr \left( {1 \over 2 \pi} \int_0^{2 \pi} \left| h_W\left(r e^{i \theta} \right) \ d \theta \right| \ \geq \ 2 \tau \right) \ \leq \  {1 \over \tau}.$$
Consequently, from (6.1) and (6.2), we have 
$$\Pr \left( \sum_{s=1}^N \ln {|a_{s, W}| \over r} \ \leq \ -{1 \over 2} \ln 2 \tau \right) \ \leq \ {1 \over \tau}.$$
and the proof follows.
{\hfill \hfill \hfill} \qed

\enddocument

\end